\documentclass[aps,11pt]{revtex4}
\usepackage{graphicx}
\usepackage{color}
\usepackage{bm}
\usepackage{amssymb}
\begin{document}

\title{Adapting to a Changing Environment: Non-Obvious Thresholds in Multi-Scale Systems}


\author{Clare Perryman and Sebastian Wieczorek}

\affiliation{
Mathematics Research Institute, University of Exeter, EX4 4QF, UK
}

\begin{abstract}
Many natural and technological systems fail to adapt to changing
external conditions and move to a different state if the conditions vary too
fast. Such ``non-adiabatic'' processes are ubiquitous, but little
understood. We identify these processes with a new nonlinear
phenomenon---an intricate {\it threshold} where a forced system fails
to adiabatically follow a changing stable state.  In systems with
multiple timescales, we derive existence conditions that show such
thresholds to be generic, but non-obvious,
meaning they cannot be captured by traditional stability theory.
Rather, the phenomenon can be analysed using concepts from modern singular
perturbation theory: folded singularities and canard trajectories,
including {\it composite canards}.  Thus, non-obvious thresholds
should explain the failure to adapt to a changing environment in a
wide range of multi-scale systems including: tipping points in the
climate system, regime shifts in ecosystems, excitability in nerve
cells, adaptation failure in regulatory genes, and adiabatic switching in technology.

\begin{description}
\item[Keywords]rate-induced bifurcations, canards, folded singularity, thresholds
\end{description}
\end{abstract}
\maketitle

\section{Introduction}

The time evolution of real-world systems often takes place on multiple
timescales, and is paced by aperiodically changing
external conditions. Of particular interest are situations where, if
the external conditions change too fast, the system fails to adapt and
moves to a different state. In climate science and ecology one speaks of
``rate-induced tipping
points''~\cite{Wieczorek2010,Lenton2008,Stocker1997,Leemans2004}, the
``critical rate hypothesis''~\cite{Scheffer2008}, and ``adaptation
failure''~\cite{Bridle2007} to describe the sudden transitions caused by
too rapid changes in external conditions (e.g. dry and hot climate
anomalies or wet periods due to El Ni\~no-Southern Oscillation).  In
neuroscience, type III excitable nerves~\cite[Ch. 7]{Izhikevich2007}
accommodate slow changes in an externally applied voltage, but
an excitation requires a rapid enough increase in the
 voltage~\cite{Hill1936,Biktashev2002}.
 In nonequilibrium genetic circuits,  cells are forced to decide 
 between alternative fates in response to changing extracellular conditions, 
 and the decision is determined by the rate of change  \cite{Nene2012}.
However, such rate-induced
transitions cannot, in general, be explained by traditional stability theory, and
require an alternative approach.

\begin{figure}[t]
\centering
\includegraphics[width=13cm]{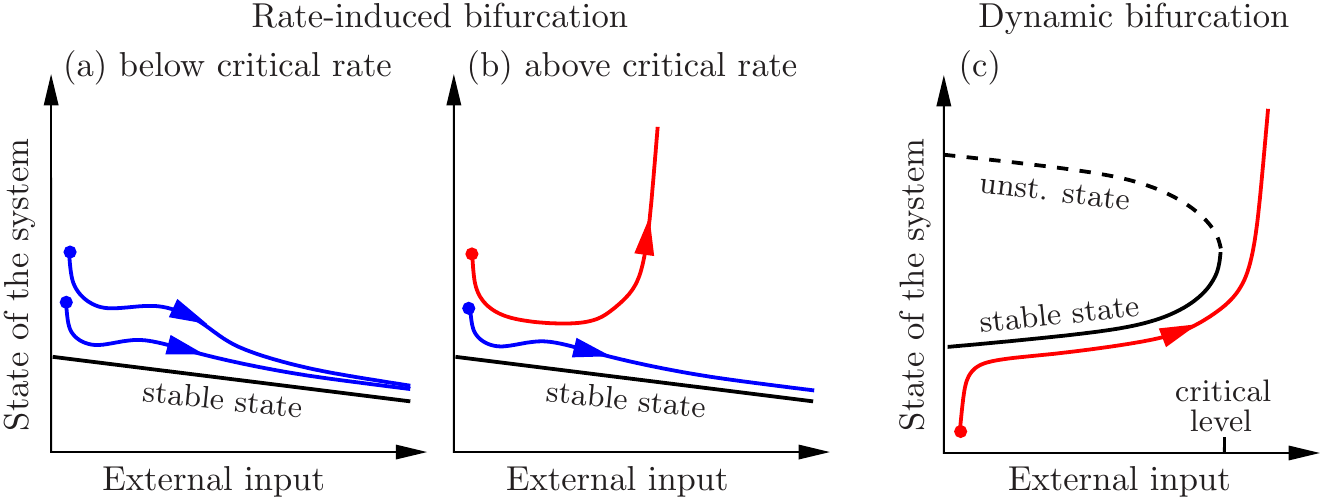}
\caption{ (Colour online) The conceptual difference between (a)--(b) a
rate-induced bifurcation and (c) a dynamic bifurcation in systems with
a time-varying external input. The ``stable state'' is an asymptotic
stable state when the external input is fixed in time.
In (a)--(b), a response to a varying external input is (red)
non-adiabatic, meaning the system destabilises, or (blue) adiabatic,
meaning the system tracks the moving stable state.  }
\label{fig:RBtip}
\end{figure}

This paper conceptualises the failure to adapt to a changing
environment as a {\it rate-induced
bifurcation}~\cite{Wieczorek2010,Ashwin2012}---a non-autonomous
instability characterised by {\it critical rates} of external
forcing~\cite{Wieczorek2010,Ashwin2012} and {\it instability
thresholds}~\cite{Wieczorek2010,Mitry2013}.  Rate-induced bifurcations
can be counter-intuitive because they occur in systems where a stable
state exists continuously for all fixed values of the external input
[Fig.~\ref{fig:RBtip}(a)--(b)].  When the external input varies in
time, the position of the stable state changes and the system tries to
keep pace with the changes. The forced system {\it adiabatically
follows} or {\it tracks} the continuously changing stable state if the
external input varies slowly enough [Fig.~\ref{fig:RBtip}(a)].
However, many systems fail to track the changing stable state if the
external input varies too fast. These systems have initial states that
destabilise---move away to a different, distant state---above some critical
rate of forcing [Fig.~\ref{fig:RBtip}(b)].  This happens even though
there is no obvious loss of stability.  Moreover, in systems with
multiple timescales there may be no obvious threshold separating the
adiabatic and non-adiabatic responses in Fig.~\ref{fig:RBtip}(b).
This is in contrast to {dynamic bifurcations}~\cite{Lobry1991}, which
can be explained by classical bifurcations of the stable state at some
{critical level} of external input [Fig.~\ref{fig:RBtip}(c)]. In this
case, the forced system destabilises totally predictably around the
critical level, independently of the initial state and of the rate of
change.

In the absence of an obvious threshold, scientists are often puzzled
by the actual boundary separating initial states that adapt to
changing external conditions from those that fail to adapt.  The first
non-obvious threshold was identified only recently, in the context of
a rate-induced climate tipping point termed the ``compost-bomb
instability'', as a {\it folded saddle canard}~\cite{Wieczorek2010}.
This finding explained a sudden release of soil carbon from peat lands
into the atmosphere above some critical rate of warming, which puzzled
climate carbon-cycle scientists~\cite{Luke2011,Wieczorek2010}.
Subsequently, similar non-obvious ``firing thresholds'' explained the
spiking behaviour of type III
neurons~\cite{Mitry2013,Wechselberger2014}.

Here, we reveal a non-obvious threshold with an intricate band
structure, and discuss the underlying mathematical mechanism.  The uncovered threshold is generic, and should explain the
failure to adapt to a changing environment in a wide range of
nonlinear multi-scale systems. Specifically, the intricate band
structure is shown to arise from a combination of the complicated
dynamics due to a folded node singularity~\cite{Szmolyan2001} and the
simple threshold behaviour due to a folded saddle
singularity~\cite{Wieczorek2010} near a folded
saddle-node type I~\cite{Krupa2010,Guckenheimer2008,Vo2014}.  What is more,
the threshold is identified with special {\it composite
canards}---trajectories that follow canard segments of different
folded singularities. More generally, we derive existence results for
critical rates and non-obvious thresholds, and discuss our
contribution in the context of canard theory and its applications.

\vspace{5mm}
\section {A general framework and existence results for non-obvious thresholds}
\label{sec:2}

Our general framework is based on geometric singular perturbation
theory~\cite{Fenichel1979,Jones1995}. It builds on the ideas developed
in~\cite{Wieczorek2010}, and extends the analysis to any type of 
smoothly varying external input. Specifically, we consider multi-scale dynamical systems akin to
simple climate, neuron, and electrical circuit
models~\cite{Luke2011,Wieczorek2010,Roberts2013,Cessi1994,Mitry2013,Wechselberger2014,Pol1934}:
\begin{eqnarray}
\delta\;dx/dt &=& f(x,y,\lambda(\epsilon t),\delta),
\label{eq:ofast}\\
dy/dt &=& g(x,y,\lambda(\epsilon t),\delta),
\label{eq:oslow}
\end{eqnarray}
with a fast variable $x$, slow variable $y$, and sufficiently smooth
functions $f$ and $g$.  The small parameter $0<\delta \ll 1$
quantifies the ratio of the $x$ and $y$ timescales.  The time-varying
external input $\lambda(\epsilon t)$ is bounded between
$\lambda_{\min}$ and $\lambda_{\max}$, and evolves smoothly on a slow timescale
$$
\tau=\epsilon t,
$$
where $\tau\in(\tau_{\min},\tau_{\max})$ can be unbounded.

The system has two small
parameters: $\delta$ and $\epsilon$.  While the analysis of
rate-induced bifurcations is greatly facilitated by the singular limit
$\delta=0$, it requires nonzero $\epsilon$. The limit $\epsilon=0$
gives the conceptual starting point for the analysis.

\begin{figure*}[t]
\centering
\includegraphics[width=13.5cm]{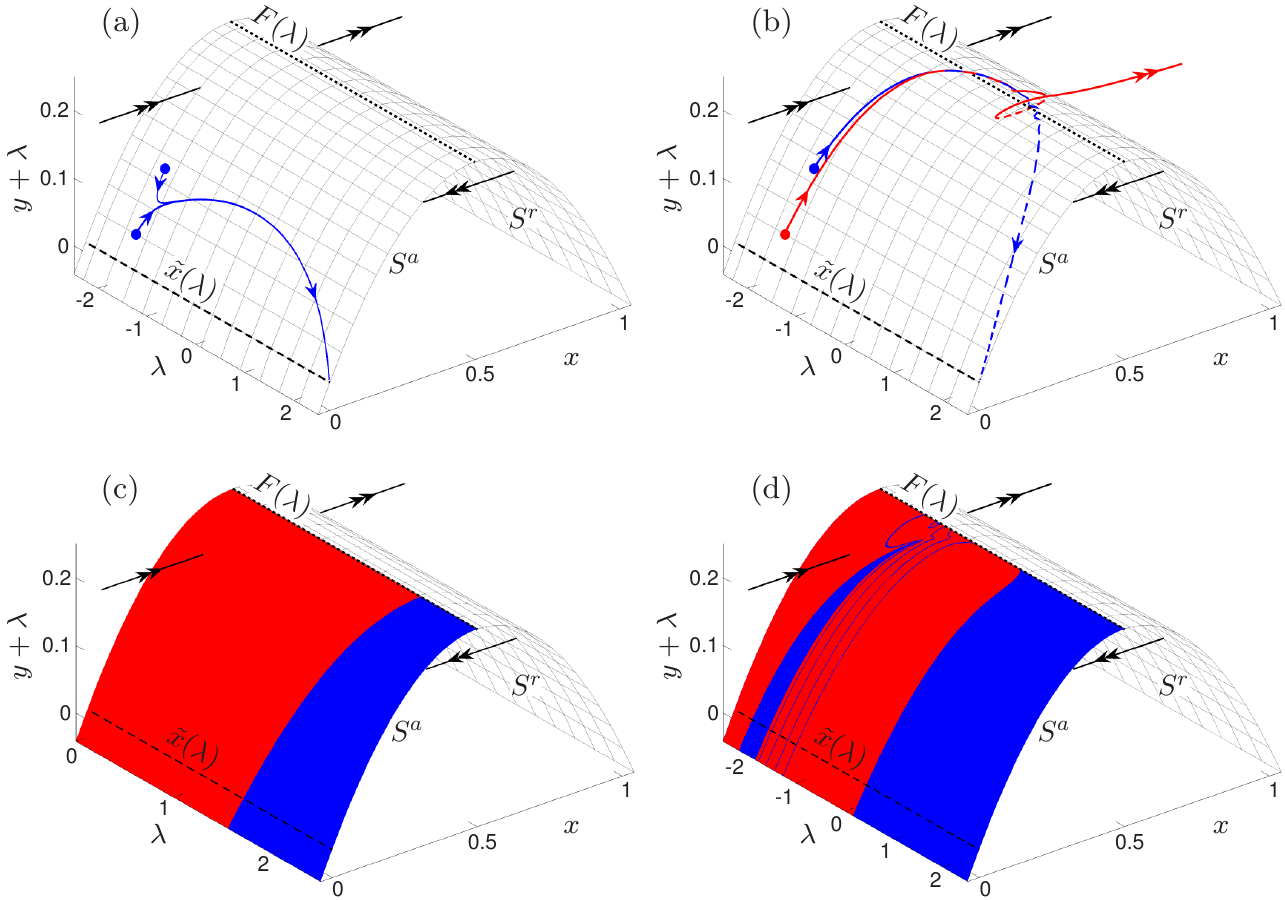}
\caption{ (Colour online) (a)--(b) Trajectories starting at two
different initial states (dots) on $S^a$, near the changing stable
state $\tilde{x}$, (a) below and (b) above the critical rate.
(c)--(d) Above the critical rate, the initial states on $S^a$ that (red)
destabilise or (blue) track $\tilde{x}(\lambda(\tau))$ highlight
different threshold types.  We used
Eqs.~(\ref{eq:ofast})--(\ref{eq:oslow}),~(\ref{eq:example}), and [(a),
(b), (d)] Eq. (\ref{eq:logistic}) with (a) $\epsilon=0.06$ and [(b),
(d)] $\epsilon=0.216$; and (c) Eq.~(\ref{eq:exponential}) with
$\epsilon=1$.  Other parameters were $\delta=0.01$,
$\lambda_{\max}=2.5$.  For (a)--(d)
the critical manifold $S(\lambda)$ is given by $y=-\lambda-x(x-1)$,
has a fold $F(\lambda)$ at $(x,y)=(1/2,-\lambda+1/4)$ and a unique
stable steady state $\tilde{x}(\lambda)$ at $(x,y)=(0,-\lambda)$.  For
clarity, the plots are shown in the co-moving coordinate system
$(x,y+\lambda,\lambda)$.  The $\lambda$ axis can be transformed into a
slow time axis using [(a), (b), (d)] Eq.~(\ref{eq:logistic}) or (c)
Eq.~(\ref{eq:exponential}). }
\label{fig:Rtiplg}
\end{figure*}

When $\lambda$ does not vary in time, i.e. when $\epsilon=0$,
Eqs.~(\ref{eq:ofast})--(\ref{eq:oslow}) define a dynamical system with
one fast and one slow variable, and a parameter $\lambda$. In the
singular limit $\delta=0$, the slow subsystem $dy/dt =
g(x,y,\lambda,0)$ evolves on the one-dimensional critical manifold
$S(\lambda)$, defined by $f(x,y,\lambda,0)=0$.  Alternatively,
$S(\lambda)$ consists of steady states of the fast subsystem $dx/dT =
f(x,y,\lambda,0)$, where $T=t/\delta$ is the fast timescale, and $y$
acts as a second parameter. The critical manifold can have an
attracting part $S^a(\lambda)$ and a repelling part $S^r(\lambda)$,
which are separated by a fold point $F(\lambda)$ (Fig.~\ref{fig:Rtiplg}).
To give precise statements
about non-obvious thresholds we assume for every fixed $\lambda$
between $\lambda_{\min}$ and $\lambda_{\max}$:
\vspace{2mm}\\
{ {\bf(A1)} \it The system has a quadratic nonlinearity.}  The
critical manifold $S(\lambda)$ is locally a graph over $x$ with a
single fold $F(\lambda)$ tangent to the fast $x$-direction, defined by
\begin{equation}
\frac{\partial f}{\partial x}\Big|_S = 0 \quad \mbox{and} \quad
\frac{\partial^2 f}{\partial x^2}\Big|_S \ne 0.
\label{eq:fold}
\end{equation}
\vspace{2mm}\\
 {{\bf(A2)} \it The system has a stable state for all fixed external
conditions.} Near $F(\lambda)$, $S^a(\lambda)$ contains just one
steady state $\tilde{x}(\lambda)$ which is asymptotically stable and
varies continuously with $\lambda$.
\vspace{2mm}\\
The geometric structure of the phase space in the singular limit
$\delta=0$ gives insight into the dynamics for $\delta$ small, but
nonzero.  Specifically, where steady states of the fast
subsystem are hyperbolic (i.e. on $S^a(\lambda)$ and $S^r(\lambda)$
but not on $F$), system~(\ref{eq:ofast})--(\ref{eq:oslow}) with
$0<\delta \ll 1$ has a slow attracting manifold
$S^a_{\delta}(\lambda)$ and a slow repelling manifold
$S^r_{\delta}(\lambda)$. Both $S^a_{\delta}(\lambda)$ and
$S^r_{\delta}(\lambda)$ are locally invariant, lie close to, and have
the same stability type as $S^a(\lambda)$ and $S^r(\lambda)$,
respectively.  This follows from Fenichel's
Theorem~\cite{Fenichel1979,Jones1995}.

When $\lambda$ varies smoothly in time such that $0<\epsilon\lesssim 1$
and $0<\delta\ll 1$,
Eqs.~(\ref{eq:ofast})--(\ref{eq:oslow}) define a dynamical system with
one fast and two slow variables:
\begin{eqnarray}
\delta\epsilon\;dx/d\tau &=& f(x,y,\lambda(\tau),\delta),
\label{eq:fast}\\
\epsilon\;dy/d\tau &=& g(x,y,\lambda(\tau),\delta),
\label{eq:slow}\\
d\tau/d\tau &=& 1.
\label{eq:t}
\end{eqnarray}
Now the critical manifolds $S^a$ and $S^r$, as well as the slow
manifolds $S^a_{\delta}$ and $S^r_{\delta}$ are two-dimensional, and
$\tilde{x}$ and $F$ form curves (Fig.~\ref{fig:Rtiplg}).  When
$\lambda(\tau)$ varies slowly enough, the forced
system~(\ref{eq:ofast})--(\ref{eq:oslow}) tracks the continuously
changing stable state $\tilde{x}(\lambda(\tau))$.  However, the
system may fail to track, and destabilise. To be more precise, we
define:
\vspace{2mm}\\
{\bf Definition 1}. For a given initial state on $S^a_{\delta}$, we
say that system~(\ref{eq:ofast})--(\ref{eq:oslow}) {\it destabilises}
if the trajectory leaves $S^a_{\delta}$ and moves away along the fast
$x$-direction.  Otherwise, we say that
system~(\ref{eq:ofast})--(\ref{eq:oslow}) {\it tracks} the moving
stable state $\tilde{x}(\lambda(\tau))$.
\vspace{2mm}\\
%
{\bf Definition 2}. The {\it critical rate} $\epsilon_c$ is the
largest $\epsilon$ below which there are no initial states on
$S^a_\delta$ that destabilise.
\vspace{2mm}\\
{\bf Definition 3}. The {\it instability threshold} is the boundary
within $S^a_{\delta}$ separating initial states that track
$\tilde{x}(\lambda(\tau))$ from those that destabilise.
\vspace{2mm}\\
Figure~\ref{fig:Rtiplg} (a)--(b) shows two trajectories of
Eqs.~(\ref{eq:ofast})--(\ref{eq:oslow}) for different initial states
on $S^a$. Below the critical rate, all trajectories track, and
eventually converge to $\tilde{x}(\lambda(\tau))$
[Fig.~\ref{fig:Rtiplg}(a)]. However, above the critical rate there are
initial states near $\tilde{x}$ that fail to track
$\tilde{x}(\lambda(\tau))$, and the system destabilises [red in
Fig.~\ref{fig:Rtiplg}(b)].  Interestingly, some trajectories leave
$S^a_{\delta}$ but, instead of destabilising along the fast
$x$-direction, return to $S^a_{\delta}$ and converge to $\tilde{x}$
[blue in Fig.~\ref{fig:Rtiplg}(b)].  The two qualitatively different
behaviours in Fig.~\ref{fig:Rtiplg}(b) show there is an instability
threshold within $S^a_{\delta}$.  What is more, the threshold can be
simple [Fig.~\ref{fig:Rtiplg}(c)] as reported
in~\cite{Wieczorek2010,Mitry2013}, or can have an intriguing band
structure [Fig.~\ref{fig:Rtiplg}(d)] that has not been reported to
date.  In both cases, it is not immediately obvious what determines
the threshold.

%
The analysis of the mathematical mechanism for non-obvious thresholds is
greatly facilitated by the singular limit $\delta=0$, where the fold and slow
manifolds are unique and known exactly.
System~(\ref{eq:fast})--(\ref{eq:t}) is reduced to the slow dynamics
on $S$ by setting $\delta=0$, and then projected onto the
$(x,\tau)$-plane by differentiating Eq.~(\ref{eq:fast}) with respect
to slow time $\tau$:
\begin{eqnarray}
 dx/d\tau &=& -\frac{g\,\partial f/\partial y + \epsilon(\partial f/\partial \lambda)(d \lambda/d\tau)}{\epsilon\,\partial f/\partial x}\Big|_{S},
\label{eq:fastp}\\
d\tau/d\tau&=& 1.
\label{eq:tp}
\end{eqnarray}
It now becomes clear that if a trajectory deviates too much from
$\tilde{x}$ and approaches a typical point on $F$ then,
according to fold condition~(\ref{eq:fold}), $\partial f/\partial x$
in Eq.~(\ref{eq:fastp}) approaches zero, and $x$ diverges off to
infinity in finite slow time $\tau$.  However, there may be special
points on $F$ where
\begin{equation}
\left[g\,\partial f/\partial y + \epsilon(\partial f/\partial \lambda)(d\lambda/d\tau) \right]|_F =0,
\label{eq:foldsing}
\end{equation}
and $dx/d\tau$ remains finite.  Such special points are referred to as
{\it folded singularities}~\cite{Takens1976,Szmolyan2001}. The
corresponding trajectories, that cross from $S^a$ along the
eigendirections of a folded singularity onto $S^r$, are referred to as
{\it singular canards}~\cite{Szmolyan2001}.  The
distinction between systems that have a critical rate and
those that do not appears to be whether there are trajectories started
on $S^a$ that reach $F$ away from a folded singularity, or whether all
trajectories started on $F$ flow onto $S^a$.  Furthermore, canard
trajectories, being solutions that separate these two behaviours, are
candidates for non-obvious thresholds.

An obstacle to the analysis of critical rates and instability
thresholds is that the flow on $F$, specifically the right hand side
of Eq.~(\ref{eq:fastp}), is not well defined. This obstacle can be
overcome by a special time rescaling~\cite{Dumortier1996}:
$$
d\tau=-ds\;\epsilon\,\left.(\partial f/\partial x)\right|_{S},
$$
where the new time $s$ passes infinitely faster on $F$, and reverses
direction on $S^r$:
  \[  
  \lim_{(x,\tau) \to F} \frac{ds}{d\tau} = \left\{
  \begin{array}{r l}
    \infty &  \text{ if } (x,\tau) \in S^a,\\
    -\infty &  \text{ if } (x,\tau) \in S^r.
  \end{array} \right.
  \]
 This gives the desingularised system
\begin{eqnarray}
dx/ds &=& \left.[g\,\partial f/\partial y + \epsilon\,(\partial f/\partial \lambda)(d \lambda/d\tau)]\right|_{S},
\label{eq:fastd}\\
d\tau/ds&=& -\epsilon\,\left.(\partial f/\partial x)\right|_{S},
\label{eq:td}
\end{eqnarray}
where trajectories remain the same as in
Eqs.~(\ref{eq:fastp})--(\ref{eq:tp}), the vector field on $F$ becomes
well defined, folded singularities become regular steady states, and
singular canards become trajectories tangent to an eigenspace of a
steady state.  One speaks of ``folded nodes'', ``folded saddles'' and
``folded foci'' for Eqs.~(\ref{eq:fastp})--(\ref{eq:tp}) if a steady
state for Eqs.~(\ref{eq:fastd})--(\ref{eq:td}) has real eigenvalues
with the same sign, real eigenvalues with opposite signs, and complex
eigenvalues with nonzero real parts, respectively.  Most importantly,
the difference between tracking and destabilising can easily be
analysed using Eqs.~(\ref{eq:fastd})--(\ref{eq:td}). Specifically, we
derive conditions for the existence of critical rates and non-obvious
thresholds:
\vspace{2mm}\\
{\bf Theorem 1.} {\it Existence of critical rates: a dissipative
Adiabatic Theorem}.  Suppose the forced
system~(\ref{eq:ofast})--(\ref{eq:oslow}) with assumptions (A1)--(A2)
satisfies the folded singularity condition~(\ref{eq:foldsing}) for
some $\tau\in(\tau_{\min},\tau_{\max})$ and
$\epsilon>0$.  Then, system~(\ref{eq:ofast})--(\ref{eq:oslow}) has a
critical rate $\epsilon_c$. The critical rate is approximately the
largest $\epsilon$ below which~(\ref{eq:foldsing}) is never satisfied
within $(\tau_{\min},\tau_{\max})$:
$$
\epsilon_c \approx \inf\left\{\epsilon>0: \left[g\,\partial f/\partial y + \epsilon(\partial f/\partial \lambda)(d\lambda/d\tau) \right]|_{F} = 0\right\}.
$$
\vspace{2mm}\\
 {\bf Theorem 2.} {\it Existence of non-obvious
thresholds}.  The forced system (\ref{eq:ofast})--(\ref{eq:oslow})
with assumptions (A1)--(A2) is guaranteed to have an
instability threshold if a folded saddle is the only folded
singularity within $(\tau_{\min},\tau_{\max})$. Then, the threshold is
given by the folded saddle maximal canard.
If $\tau_{\max} = \infty$ and $\lambda(\tau)$ is asymptotically constant
\begin{equation}
 \lim_{\tau \to \infty} \frac{d\lambda}{d\tau}=0,
\label{eq:asymptotic}
\end{equation}
then the system has an instability threshold if, and only if, there is
a folded saddle singularity.
\vspace{2mm}\\
{\bf Note.}  Often in real-life applications the changing external
conditions $\lambda$ are expressed as a prescribed function of time
$t$, but not $\epsilon$ or $\tau$.  Specifying $\epsilon$ is not
necessary. If one replaces $\tau$ with $\epsilon t$ in
Eqs.~(\ref{eq:fast})--(\ref{eq:td}) the dependence on $\epsilon$
disappears. However, $\epsilon$ and $\tau$ are useful for defining
critical rates of change, and facilitate the derivation of the
statements in Theorems 1 and 2.
\vspace{2mm}\\

The proofs, given in the Appendix, are based on two steps.  In the
first step, a qualitative analysis of
Eqs.~(\ref{eq:fastd})--(\ref{eq:td}) identifies the appearance of a folded singularity
with a critical rate, and certain singular canards as candidates for an
instability threshold.  In the second step, recent
results from canard
theory~\cite{Szmolyan2001,Wechselberger2005,Vo2014} are used that
state singular canards due to folded saddles, folded nodes, and
folded saddle-nodes of type I, perturb to {\it maximal canards}
in~(\ref{eq:fast})--(\ref{eq:t}) with $0<\delta \ll 1$.  Maximal
canards are those trajectories crossing from $S^a_{\delta}$ onto
$S^r_{\delta}$, which remain on $S^r_{\delta}$ for the longest time.
In this paper, we numerically compute both maximal canards
$\gamma_{\delta}$, shown in Fig.~\ref{fig:case2b}, and their
approximations by singular canards $\gamma$, shown in
Figs.~\ref{fig:case2a} and \ref{fig:case1}.

\vspace{5mm}
\section{Two cases of a non-obvious threshold}

Guided by the proof of Theorem 2, specifically the analysis of the
phase portraits containing a folded saddle [Appendix,
Fig.~\ref{fig:phaseportraits}(a)--(b)], we distinguish two cases of a
non-obvious threshold.  Furthermore, we identify one case with the
complicated threshold shown in Fig.~\ref{fig:Rtiplg}(d), and uncover
the underlying mechanism.

We illustrate the two cases using an example
of~(\ref{eq:ofast})--(\ref{eq:oslow}) with
\begin{eqnarray}
f=x(x-1)+y+\lambda(\tau)\;\;\;\mbox{and}\;\;\;g=-x,
\label{eq:example}
\end{eqnarray}
and two different aperiodic forcing functions $\lambda(\tau)$
satisfying~(\ref{eq:asymptotic}).

{\it Case 1: Complicated threshold due to a folded saddle-node type I
singularity.}  Consider example~(\ref{eq:example}) subject to logistic
growth at a rate $\epsilon$:
\begin{eqnarray}
\lambda(\tau)=\lambda_{\max}\,\tanh\!\left(\tau\right),
\label{eq:logistic}
\end{eqnarray}
where $\lambda\in(-\lambda_{\max},\lambda_{\max})$, $\tau\in(-\infty,\infty)$
and $\tau=\epsilon t$.  The desingularised
system~(\ref{eq:fastd})--(\ref{eq:td}) becomes
\begin{eqnarray}
dx/ds &=& -x + \frac{\epsilon}{\lambda_{\max}}\left(\lambda_{\max}^2 -
\lambda^2(\tau)\right), \label{eq:logisticdfast}\\
d\tau/ds &=& \epsilon(1 - 2x).
\label{eq:logisticd}
\end{eqnarray}
Steady states of~(\ref{eq:logisticdfast})--(\ref{eq:logisticd}) lie on the fold $x=1/2$, at $\lambda(\tau)$
satisfying the folded singularity condition~(\ref{eq:foldsing}):
\begin{equation}
\lambda^2(\tau) - \lambda_{\max}\left(\lambda_{\max} -\frac{1}{2\epsilon} \right)=0,
\label{eq:foldsing1}
\end{equation}
and their eigenvalues $\xi$ are found from the characteristic polynomial
\begin{equation}
\xi^2 + \xi - 4\epsilon^2\lambda(\tau)\left[
1-\left(\frac{\lambda(\tau)}{\lambda_{\max}}\right)^2\right]=0.
\label{eq:foldsing1b}
\end{equation}
The folded singularity condition~(\ref{eq:foldsing1}) has no real
roots when $\epsilon < (2\lambda_{\max})^{-1}$. 
When $\epsilon = (2\lambda_{\max})^{-1}$, there is a double root
within $(\tau_{\min},\tau_{\max})$, corresponding to a folded
saddle-node type I~\cite{Krupa2010} at $(x,\lambda(\tau))=(1/2,\,0)$.  When
$\epsilon > (2\lambda_{\max})^{-1}$, there are two distinct roots
within $(\tau_{\min},\tau_{\max})$, corresponding to a stable folded
node (focus) $FN(FF)$ at
$(x,\lambda(\tau))=(1/2,\,-\sqrt{\lambda_{\max}(\lambda_{\max}-(2\epsilon)^{-1})}\,)$
and a folded saddle $FS$ at
$(x,\lambda(\tau))=(1/2,\,\sqrt{\lambda_{\max}(\lambda_{\max}-(2\epsilon)^{-1})}\,)$.
This means that, upon increasing $\epsilon$, there is a generic saddle
node bifurcation of folded singularities at
$\epsilon_{SN}=(2\lambda_{\max})^{-1}$, which by Theorem 1 is
approximately the critical rate $\epsilon_c$ for $0<\delta\ll 1$.
According to Theorem 2, condition~(\ref{eq:asymptotic}) and the
presence of a folded saddle guarantee an instability threshold.
However, unlike the case of an isolated folded saddle whose threshold is specified by Theorem 2, it is not immediately
clear what forms the threshold near a folded saddle-node type I~\cite{Vo2014}.  
Nonetheless, this can be established numerically.

\begin{figure*}[t]
\centering
\includegraphics[width=11cm]{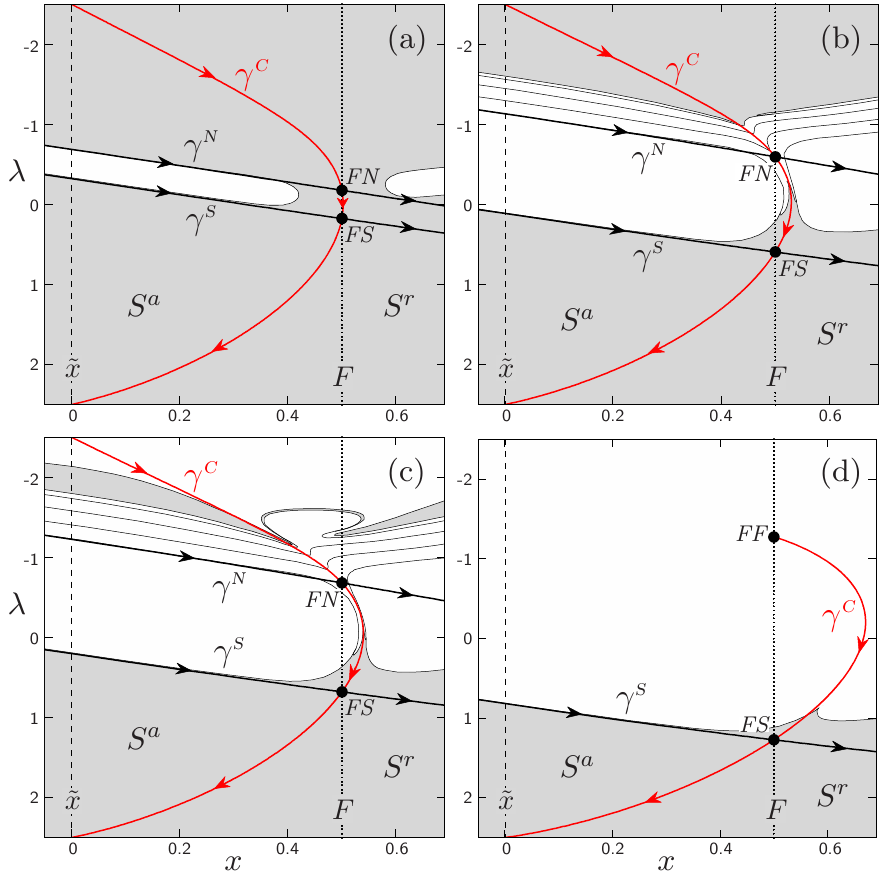}
 \caption{
(Colour online) Initial states on the critical manifold $S$ that
(white) destabilise or (grey) track $\tilde{x}(\lambda(\tau))$
in Eqs.~(\ref{eq:ofast})--(\ref{eq:oslow}) and
(\ref{eq:example})--(\ref{eq:logistic}) with $\delta=0.01$,
$\lambda_{\max}=2.5$, and $\epsilon=
\mbox{(a)}\,0.201,\,\mbox{(b)}\,0.212,\,\mbox{(c)}\,0.216,\,\mbox{and
(d)}\,0.270$, shown projected onto the $(x,\lambda)$ plane.
 Away from $F$, the instability threshold in
$S^a_{\delta}$ is well approximated by the white-grey boundary in
$S^a$.  Points $FN$, $FS$, and $FF$ are folded node, folded saddle,
and folded focus singularities, respectively; the strong folded node
singular canard $\gamma^N$ and the folded saddle singular canard
$\gamma^S$ approximate projections of the maximal canards
$\gamma^N_{\delta}$ and $\gamma^S_{\delta}$, respectively, onto $S$.
The projection of the maximal canard $\gamma^C_{\delta}$ onto  $S$ is
approximated by the weak folded node/faux saddle singular canard $\gamma^C$
when $\lambda > -1$, but lies below $\gamma^C$ for $-2.5<\lambda<-1$, e.g.
within the wide grey band around $\lambda=-2$ in (c).
Although it is difficult to see, $\gamma^C$
terminates on $F$ just above $FF$.
Compare (c) with
Fig.~\ref{fig:Rtiplg}(d).
 }
\label{fig:case2a}
\end{figure*}

The instability threshold is defined on the attracting slow manifold
$S^a_\delta$, which is difficult to compute near the fold
$F$.  To facilitate numerical computations, we consider
initial states on the critical manifold $S^a$, which is known exactly.
The results are shown in Fig.~\ref{fig:case2a}, where the white
regions indicate destabilising, and the grey regions indicate
tracking.  Away from $F$, the critical manifold $S^a$ closely
approximates the slow manifold $S^a_\delta$. Here, the instability
threshold is well approximated by the boundaries between the white and
grey regions. However, caution is required near $F$, especially
around $FN$, where $S^a_\delta$ twists in a complicated
manner~\cite[Fig. 6]{Desroches2012}, and the chosen surface of initial
conditions, $S^a$, intersects these twists.  There, the boundaries
between the white and grey regions deviate from the instability
threshold due to the choice of initial states.  We also show what
happens to initial states on $S^r$ just to the right of $F$, as some
are mapped along the fast flow onto $S^a_\delta$ and converge to
$\tilde{x}$. This is why a ``reflection'' of the band structure from
$S^a$ can be seen on $S^r$.

%

Shortly past the saddle-node bifurcation, there are three bands of
initial states on $S^a$ [Fig.~\ref{fig:case2a}(a)].  The
threshold separating these bands is formed by two canard trajectories:
the folded saddle maximal canard $\gamma^S_{\delta}$, and the strong
folded node maximal canard $\gamma^N_{\delta}$. On $S^a$, trajectories
started in the white band enclosed by $\gamma^S_{\delta}$ and
$\gamma^N_{\delta}$ move directly towards the fold, then leave the
attracting slow manifold $S^a_{\delta}$ and destabilise along the fast
$x$-direction.  Trajectories started in the grey band below
$\gamma^S_\delta$ approach the faux saddle maximal canard
$\gamma^C_\delta$ straight away, thereby staying on the attracting
slow manifold $S^a_\delta$ and tracking $\tilde{x}$.  This is in
contrast to trajectories started in the other grey band on $S^a$, the one
above $\gamma^N_\delta$.  These trajectories initially approach and
twist around the weak folded node maximal canard $\gamma^C_\delta$,
and leave $S^a_{\delta}$. However, rather than destabilising, they are
fed back along $\gamma^C_\delta$, onto $S^a_{\delta}$, and eventually
remain on $S^a_\delta$ [Fig.~\ref{fig:Rtiplg}(b), blue trajectory].
Finally, grey initial states on $S^r$ are mapped along the
fast flow onto the grey bands on $S^a_\delta$.

\begin{figure}[t]
\centering
\includegraphics[width=13.5cm]{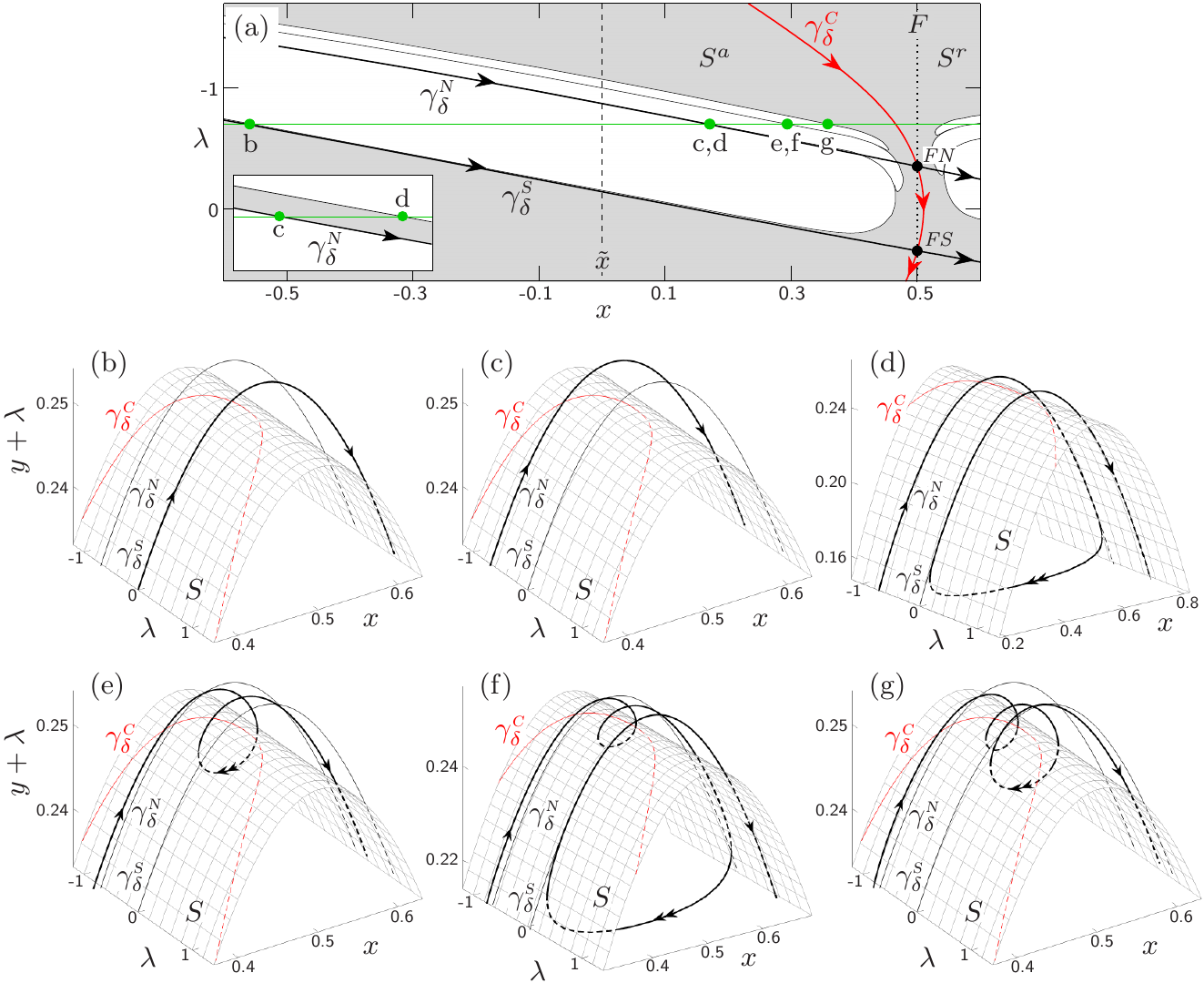}
\caption{ (Colour online) (a) Initial states on the critical manifold
$S$ that (white) destabilise or (grey) track
$\tilde{x}(\lambda(\tau))$ in
Eqs.~(\ref{eq:ofast})--(\ref{eq:oslow})
and~(\ref{eq:example})--(\ref{eq:logistic}) with $\delta=0.01$ and
$\epsilon=0.204$. Inset shows grey band between c and d; a similar
band exists between e and f.  Labels b--g at $\lambda=-0.7$, or $\tau=-\tanh^{-1}(0.28)$,
denote different threshold components including: (b) the folded saddle
maximal canard $\gamma^S_{\delta}$, (c) the strong folded node maximal
canard $\gamma^N_{\delta}$, (d) a composite canard that follows
$\gamma^N_{\delta}$ and $\gamma^S_{\delta}$, (e) a secondary folded
node maximal canard, (f) a composite canard that follows a secondary
maximal canard and $\gamma^S_{\delta}$, (g) a secondary folded node
maximal canard.  }
\label{fig:case2b}
\end{figure}

As $\epsilon$ increases, the threshold becomes more complicated due to
the presence of the stable folded node $FN$. Additional threshold
curves appear successively above $\gamma^N_{\delta}$, giving up to
five white bands of initial states above $\gamma^N_{\delta}$ that
destabilise [Fig.~\ref{fig:case2a}(b)]. Trajectories started within
these additional white bands twist around $\gamma^C_\delta$ before
destabilising, [Fig.~\ref{fig:Rtiplg}(b), red trajectory].  These
white bands are separated by narrow grey bands which are difficult to
see in Fig.~\ref{fig:case2a}; see the narrow grey band in the inset of
Fig.~\ref{fig:case2b}, or narrow blue bands in
Fig.~\ref{fig:Rtiplg}(d).  Trajectories started within
these narrow grey bands leave $S^a_\delta$, follow a maximal canard on
$S^r_\delta$ for some time, but then return to $S^a_\delta$ into the grey
region below $\gamma^S_\delta$, and converge to $\tilde{x}$.  The white bands expand
with $\epsilon$ and approach the weak folded node maximal canard
$\gamma^{C}_{\delta}$ on both sides [Fig.~\ref{fig:case2a}(c)].  When
the folded node $FN$ turns into a folded focus $FF$ at $\epsilon = (2
+\sqrt{4 + \lambda_{\max}^2}\,)/8\lambda_{\max}$, its canards
disappear~\cite{Szmolyan2001} and so does the band structure
[Fig.~\ref{fig:case2a}(d)]. We are left with a simple threshold, given
just by $\gamma^S_{\delta}$ as in Ref.~\cite{Wechselberger2014}.

%
The key mechanism for complicated thresholds is the phenomenon whereby
trajectories leave $S^a_\delta$ through the folded node region and
then, rather than destabilising, are fed back to $S^a_\delta$ through
the folded saddle region. This phenomenon has two
consequences. Firstly, not all initial states on $S^a_\delta$ and
above $\gamma^N_\delta$ destabilise. Secondly, the initial states on
$S^a_\delta$ that destabilise or track $\tilde{x}$ form alternating
bands, and these bands have not been identified before. More
generally, the alternating bands are related to the known
rotational sectors of a folded node; see \cite{Wechselberger2005} 
for a detailed discussion of rotational sectors. However, whilst rotational
sectors are separated by a single canard
trajectory~\cite{Desroches2012,Wechselberger2005}, our white bands are
separated by a narrow grey band bounded by two different canard
trajectories.

Figure~\ref{fig:case2b} identifies the different components of the
complicated threshold. They consist of known maximal canards such as
(b) $\gamma^S_{\delta}$, (c) $\gamma^N_{\delta}$, and [(e) and (g)]
secondary folded node maximal canards that bifurcate off
$\gamma^C_{\delta}$~\cite{Wechselberger2005}.  These canards form the
lower boundaries of the narrow grey bands. Most interestingly, they also consist of new
{\it composite canards} that follow canard segments of
different folded singularities. These canards form the upper
boundaries of the narrow grey bands. Figure~\ref{fig:case2b} shows
composite canards which initially (d) follow $\gamma^N_{\delta}$, or
(f) follow the first secondary folded node maximal canard, and then
[(d) and (f)] follow $\gamma^S_{\delta}$.  This explains the
intriguing band structure with intermingled regions of white and grey
in Figs.~\ref{fig:case2a}(b)-(c) and~\ref{fig:case2b}(a), or red and
blue in Fig.~\ref{fig:Rtiplg}(d).  The
composite canards in Fig.~\ref{fig:case2b}(d) and (f) are reminiscent
of trajectories that switch between different primary and secondary
canards of the same folded node in a stellate cell
model~\cite{Wechselberger2009} and in a reduced Hodgkin-Huxley model~\cite[Fig. 9]{Desroches2010}.

{\it Case 2: Simple threshold due to an isolated folded saddle
singularity.}  Consider example~(\ref{eq:example}) subject to an
exponential approach at a rate $\epsilon$:
\begin{eqnarray}
\lambda(\tau)=\lambda_{\max}\left(1 - e^{-\tau}\right),
\label{eq:exponential}
\end{eqnarray}
where $\lambda\in(0,\lambda_{\max})$, $\tau\in(0,\infty)$ and $\tau=\epsilon t$.  The
desingularised system~(\ref{eq:fastd})--(\ref{eq:td}) becomes
\begin{eqnarray}
dx/ds &=& -x + \epsilon\left(\lambda_{\max} - \lambda(\tau)\right), \label{eq:exponentialdfast}\\
d\tau/ds &=& \epsilon(1 - 2x).
\label{eq:exponentiald}
\end{eqnarray}
The steady state of~(\ref{eq:exponentialdfast})-(\ref{eq:exponentiald}) lies on the fold $x=1/2$, at $\lambda(\tau)$
satisfying the folded singularity condition~(\ref{eq:foldsing}):
\begin{equation}
\lambda(\tau) = \lambda_{\max} - \frac{1}{2\epsilon},
\label{eq:foldsing2}
\end{equation}
and its eigenvalues $\xi$ are found from the characteristic polynomial
\begin{equation}
\xi^2 + \xi + 2\epsilon^2\left(\lambda(\tau)-\lambda_{\max}\right)=0.
\label{eq:foldsing2b}
\end{equation}
%
%
\begin{figure}[t]
\centering
\includegraphics[width=11cm]{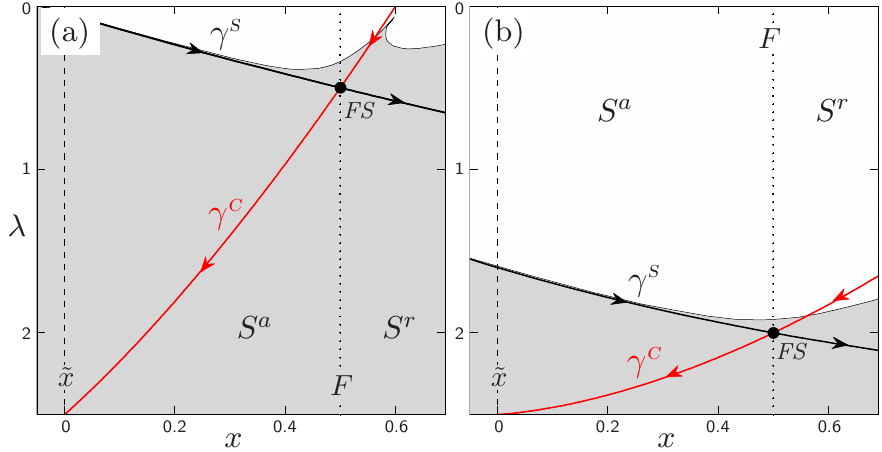}
\caption{
(Colour online) Initial states on the critical manifold $S$ that
(white) destabilise or (grey) track $\tilde{x}(\lambda(\tau))$
for Eqs.~(\ref{eq:ofast})--(\ref{eq:oslow}), (\ref{eq:example}), and
(\ref{eq:exponential}) with $\delta=0.01$, $\lambda_{\max}=2.5$, and
(a) $\epsilon=0.25$, (b) $\epsilon=1$, shown projected onto the
$(x,\lambda)$ plane. Away from $F$ the instability threshold in
$S^a_{\delta}$ is well approximated  by the white-grey boundary in
$S^a$.  Compare (b) with
Fig.~\ref{fig:Rtiplg}(c). For labels see
Fig.~\ref{fig:case2a}.}
\label{fig:case1}
\end{figure}
%
The main difference from Case 1 is that the different forcing
$\lambda(\tau)$ in (\ref{eq:exponential}) gives a folded singularity
condition~(\ref{eq:foldsing2}) with just a single root, corresponding
to an isolated folded saddle $FS$ at $(x,\lambda(\tau))
=(1/2,\,\lambda_{\max}-(2\epsilon)^{-1})$. Upon increasing $\epsilon$,
the folded saddle enters $(\tau_{\min},\tau_{\max})$ via its lower
boundary when $\epsilon=(2\lambda_{\max})^{-1}$, which by Theorem 1 is
approximately the critical rate $\epsilon_c$ for $0<\delta\ll
1$. According to Theorem 2, there is an instability threshold given by
the folded saddle maximal canard $\gamma^S_\delta$, as in the
compost-bomb and the type III neuron
examples~\cite{Wieczorek2010,Mitry2013}.  Numerical computations in
Fig.~\ref{fig:case1} confirm that for $\delta=0.01$, and away from $F$,
the threshold is well approximated by the singular canard $\gamma^S$.
It is interesting to note, the threshold in Fig.~\ref{fig:case1} is
very similar to that in Fig.~\ref{fig:case2a}(d).

{\it Note on types of non-obvious thresholds.}  
Theorem 2 in conjunction with numerical investigations in this section show that
which case of a non-obvious threshold occurs, if any at all, depends both
on the system~(\ref{eq:ofast})--(\ref{eq:oslow}) and the form of 
the external input $\lambda(\tau)$.  Specifically, the threshold is determined 
by the number, type and stability of the folded singularities.  What is more,
our simple example~(\ref{eq:example}) demonstrates that both cases of a non-obvious
threshold can occur for the same system when subject to different $\lambda(\tau)$.

In both cases, the external input $\lambda(\tau)$ 
satisfies~(\ref{eq:asymptotic}).  When $\lambda(\tau)$ does
not satisfy~(\ref{eq:asymptotic}), there can be an instability
threshold that is not associated with a folded saddle [Appendix,
Fig.~\ref{fig:phaseportraits}(d)].  However,
it follows from the proof of Theorem 2 in the Appendix that such a
threshold is simple, akin to the case of an isolated folded saddle.

\vspace{5mm}
\section{Conclusions}

In summary, we analysed multiple timescale systems subject to an
aperiodically changing environment, identified nonlinear mechanisms
for the failure to adapt, and derived conditions for the existence of
these mechanisms. Specifically, we discussed instability thresholds
where a system fails to adiabatically follow a continuously changing
stable state.  Despite their cross-disciplinary nature, these
thresholds are largely unexplored because they are ``non-obvious'',
meaning they cannot, in general, be revealed by traditional stability
theory.  Thus, they require an alternative approach.  We presented a
framework, based on geometric singular perturbation theory, that led
us to a novel type of threshold with an intriguing band structure.  The
threshold has alternating bands, where the system tracks the moving
stable state, or destabilises.  We showed that this structure is
organised by a folded saddle-node type I singularity.  Intuitively, it
arises from an interplay of the complicated dynamics of twisting
canard trajectories due to a folded node singularity, and the simple
threshold behaviour illustrated for a folded saddle singularity. Most
importantly, trajectories which leave the attracting slow manifold
through the folded node region can be fed back to the attracting slow
manifold through the folded saddle region.  In more technical terms,
the band structure is related to the rotational sectors of a folded
node, but also differs from them in one key aspect. Whereas the
rotational sectors are separated by a single canard trajectory, namely the
maximal canard~\cite{Desroches2012,Wechselberger2005}, the
corresponding wide bands are separated by a narrow band.  These
separating narrow bands are bounded by two different canard
trajectories. One of them is a known maximal canard, and the other is
a composite canard that follows maximal-canard segments of different
folded singularities.

Whilst non-obvious thresholds can be complicated, they are generic,
and should explain counter-intuitive responses to a changing
environment in a wide range of multi-scale systems. We highlighted
their importance by examples of climate and ecosystems failing to
adapt to a rapidly changing
environment~\cite{Wieczorek2010,Ashwin2012,Luke2011}, and type III
excitable cells ``firing'' only if the voltage stimulus rises fast
enough~\cite{Mitry2013,Hill1936}. More generally, our results give new
insight into non-adiabatic processes in multi-scale dissipative
systems, and should stimulate further work in canard theory.

\section*{Acknowledgements}
We would like thank M. Wechselberger and T. Vo for
useful discussions.  The research of C.P. was supported by the EPSRC
and the MCRN (via NSF grant DMS-0940363).

\section*{Appendix}
Consider system~(\ref{eq:ofast})--(\ref{eq:oslow})
with assumptions
(A1)--(A2), and restrict the discussion to
$\tau\in(\tau_{\min},\tau_{\max})$, which can be unbounded.

\subsection*{\it Proof of Theorem 1.}
Let $p$ be a point on the fold $F$ in the desingularised
system~(\ref{eq:fastd})--(\ref{eq:td}).  By (A1) the vector field at
$p$ only has a component in the $x$-direction.  When $\epsilon=0$, by
assumption (A2) the vector field points towards the attracting
critical manifold $S^a$ at every $p\in F$. This means all trajectories
starting on $F$ flow onto $S^a$, and no trajectories starting on $S^a$
reach $F$. When $\epsilon>0$, there may be trajectories that reach $F$
from $S^a$.  This happens if, and only if, the vector field changes
sign at some $p\in F$ as $\epsilon$ is varied:
\begin{eqnarray}
(dx/ds)|_{p} &=& \left.[g\,\partial f/\partial y + \epsilon\,(\partial f/\partial \lambda)(d \lambda/d\tau)]\right|_{p} = 0,
\label{eq:cross1}\\
\frac{d}{d\epsilon}(dx/ds)\Big|_{p} &=& [(\partial f/\partial \lambda)(d \lambda/d\tau)]|_{p} \ne 0.
\label{eq:cross2}
\end{eqnarray}
Furthermore, by assumption (A1) $S$ can be expressed as a graph over
$y$ meaning $(\partial f/\partial y)|_{p}\ne 0$, and by assumption
(A2) there are no steady states on $F$ in the full system meaning
$g|_{p}\ne0$, so~(\ref{eq:cross1}) already implies~(\ref{eq:cross2}).
 
By~\cite[Th. 1]{Jones1995}, if system~(\ref{eq:fastd})--(\ref{eq:td})
has no trajectories started on $S^a$ that reach $F$, then
system~(\ref{eq:fast})--(\ref{eq:t}) has no trajectories that leave
$S^a_{\delta}$ for $0<\delta\ll 1$.  Furthermore,
by~\cite[Th. 1]{Szmolyan2004}, if system~(\ref{eq:fast})--(\ref{eq:t})
has trajectories starting on $S^a$ that reach $F$ away from a folded
singularity, then system~(\ref{eq:fast})--(\ref{eq:t}) has
trajectories that leave $S^a_{\delta}$ and move away along the fast
$x$-direction for $0<\delta\ll 1$. Hence, the folded singularity
condition~(\ref{eq:cross1}) implies a critical rate for
system~(\ref{eq:fast})--(\ref{eq:t}), and for the original
system~(\ref{eq:ofast})--(\ref{eq:oslow}).

By Definition 2, in the singular limit $\delta=0$ the critical rate is the largest
$\epsilon$ below which~(\ref{eq:cross1}) is never satisfied within
$(\tau_{\min},\tau_{\max})$. When $\delta$ is small but nonzero, the critical rate is given by 
$$
\epsilon_c \approx \inf\left\{\epsilon>0: \left[g\,\partial f/\partial y + \epsilon(\partial f/\partial \lambda)(d\lambda/d\tau) \right]|_{F} = 0\right\} + E_\delta,
$$
where $E_\delta$ is a correction for nonzero $\delta$.
For $\delta$ small enough, the correction term $E_\delta$ is $O({\delta}^{\frac{1}{2}})$ if the 
folded singularity at $\epsilon_c$ is a saddle, node, or folded saddle-node type II 
\cite{Szmolyan2001, Krupa2010}, and is $O({\delta}^{\frac{1}{4}})$ if the 
folded singularity at $\epsilon_c$ is a folded-saddle node type I \cite{Vo2014}.

\noindent
\subsection*{\it Proof of Theorem 2.}
%

\begin{figure}[t]
\centering
\includegraphics[width=13.5cm]{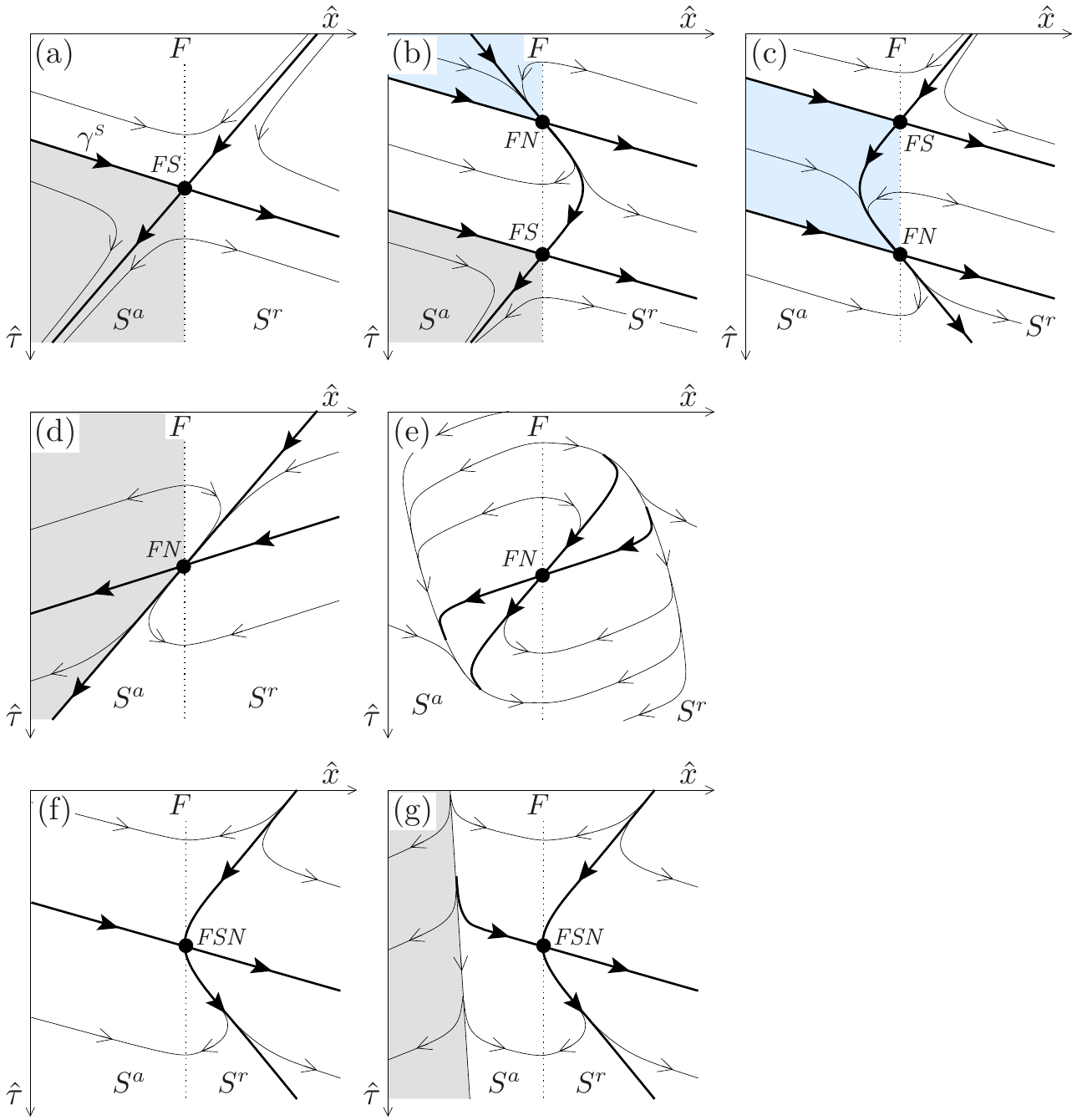}
\caption{ (Colour online) Sketches of selected phase portraits for
system~(\ref{eq:fastp})-(\ref{eq:tp}), containing folded saddles
($FS$), folded nodes ($FN$), and folded saddle-nodes ($FSN$). Singular
canards are shown in bold. On $S^a$, there are trajectories that (white)
approach $F$ away from a folded singularity, (blue) 
leave $S^a$ via a folded singularity, and (grey) never reach $F$.}
\label{fig:phaseportraits}
\end{figure}
%

Consider a fixed value of $\epsilon > \epsilon_c$.  We
are interested in phase portraits of
system~(\ref{eq:fastp})--(\ref{eq:tp}) which have two types of
trajectories starting on $S^a$: those that reach $F$ away from a
folded singularity, and those that never reach $F$ and remain on
$S^a$.  We refer to the separatrix dividing these two types of
trajectories as the {\it singular threshold}. Phase
portraits of system~(\ref{eq:fastp})--(\ref{eq:tp}) that may contain a
singular threshold are identified as follows. We keep in mind that
$d\tau/dt>0$, construct possible phase portraits of the desingularised
system~(\ref{eq:fastd})--(\ref{eq:td}), reverse the flow on $S^r$, and
keep those portraits that allow a singular threshold.

The proof consists of three parts.  Firstly, we analyse
an arbitrary external input $\lambda(\tau)$ to show that an isolated
folded saddle guarantees a singular threshold.  Secondly, we
analyse
an asymptotically constant external input, i.e. $\lambda(\tau)$ satisfies condition~(\ref{eq:asymptotic}),
to show there is a singular threshold if, and only if, there is a folded
saddle.  Lastly, we use recent results from canard theory to show that
singular thresholds persist as instability thresholds for $\delta$
small, but nonzero.

\subsubsection*{Part 1}


Firstly, assume condition (\ref{eq:foldsing}) is satisfied, meaning there
is a folded singularity $p$. Without loss of
generality, suppose $p$ is at the origin. According
to~\cite[Prop. 2.1]{Szmolyan2001}, under assumption (A1) and
condition (\ref{eq:foldsing}), there is a smooth change of
coordinates that projects the fold curve $F$ orthogonally onto the
$\tau$-axis and, in the neighbourhood of $p$, brings the
desingularised system~(\ref{eq:fastd})--(\ref{eq:td}) to the normal form

\begin{eqnarray}
   \frac{d\hat{x}}{d\hat{s}}&=&b\hat{\tau} + c\hat{x} +
   O(\hat{x}^2,\hat{x}\hat{\tau},\hat{\tau}^2) \label{eq:fastn},\\
   \frac{d\hat{\tau}}{d\hat{s}}&=&-2\epsilon\hat{x} +
   O(\hat{x}^2,\hat{x}\hat{\tau}), \label{eq:tn}
\end{eqnarray}
where $\hat{x}$ and $\hat{\tau}$ are the new coordinates, the fold $F$
is defined by $\hat{x}=0$, and the attracting critical manifold $S^a$
is defined by $\hat{x}<0$.  The eigenvalues of $p$:
$$ 
\xi_{1,2}=(c \pm \sqrt{c^2 - 8\epsilon b}\,)/2,
$$ 
determine the type of the folded singularity in
system~(\ref{eq:fastp})--(\ref{eq:tp}). In particular, $p$ is a folded
saddle if $b<0$, a folded saddle-node if $b=0$, and a folded node,
focus or centre if $b>0$. The key observation for our purposes is that
$b\ne0$ determines the direction of the flow on $F$ where
$d\hat{\tau}/d\hat{s}=0$ and $d\hat{x}/d\hat{s}= b\hat{\tau}+
O(\hat{\tau}^2)$.

In the case of a folded saddle ($b<0$), trajectories starting on $S^a$
and near $F$ reach $F$ when $-1\ll\hat{\tau}<0$, or flow away from $F$
onto $S^a$ when $0<\hat{\tau}\ll 1$
[Fig.~\ref{fig:phaseportraits}(a)].  If a folded saddle is the only
folded singularity, then there are no additional changes in the
direction of the flow on $F$. The local behaviour for $0<\hat{\tau}\ll
1$ extends to $0<\hat{\tau}<\hat{\tau}_{\max}$, meaning no
trajectories started on $S^a$ for $\hat{\tau}>0$ ever reach
$F$. Hence, an isolated folded saddle implies a singular threshold.
What is more, the threshold is given by the singular folded saddle
canard. This can be seen by noting that, in the desingularised
system~(\ref{eq:fastd})--(\ref{eq:td}), the separatrix between
trajectories starting on $S^a$ that reach $F$ and those that never
reach $F$ is the stable manifold of the saddle equilibrium.  This
stable manifold becomes the singular folded saddle canard
$\gamma_{\delta}^S$ in system~(\ref{eq:fastp})--(\ref{eq:tp})
[Fig.~\ref{fig:phaseportraits}(a)].  
If, in addition to a folded saddle, there are other folded
singularities, a singular threshold can no longer be guaranteed
 [e.g. Fig.~\ref{fig:phaseportraits}(c)], nor excluded
[e.g. Fig.~\ref{fig:phaseportraits}(b)]. To obtain the threshold, one
needs to study the behaviour of trajectories started on $S^a$;
see the analysis of Case 1 in Section~3.

In the special case of a folded saddle-node ($b=0$), the
flow on $F$ in system~(\ref{eq:fastn})--(\ref{eq:tn}) is determined by
$d\hat{x}/d\hat{s}= O(\hat{\tau}^2)$.  This means there is no change
in the sign of the flow at $p$
[e.g. Fig.~\ref{fig:phaseportraits}(f)].  A folded saddle-node is
structurally unstable. Under arbitrarily small variation of system
parameters, it unfolds into a folded saddle at positive $\hat{\tau}$
and a folded node at negative $\hat{\tau}$ (multiple singularities discussed in the
paragraph above), or into no singularities.  
In the case of a folded node, focus or centre ($b>0$), trajectories
starting on $S^a$ and sufficiently close to $F$ flow away from $F$
onto $S^a$ when $-1\ll\hat{\tau}<0$, or reach $F$ when
$0<\hat{\tau}\ll 1$; see an example of an unstable folded node in
Fig.~\ref{fig:phaseportraits}(d).  For $b \geq0$, a singular
threshold cannot be guaranteed
[e.g. Fig.~\ref{fig:phaseportraits}(f)--(g)], nor excluded
[e.g. Fig.~\ref{fig:phaseportraits}(d)--(e)].

Secondly, assume there are no folded singularities.  If
the flow on $F$ in system~(\ref{eq:fastn})--(\ref{eq:tn}) points
towards $S^a$, a singular threshold can be excluded. If the flow on
$F$ points towards $S^r$, a singular threshold cannot be guaranteed,
nor excluded [restricting the $(\hat{\tau}_{\min},\hat{\tau}_{\max})$
interval to the lower part of the phase portrait in
Fig.~\ref{fig:phaseportraits}(d) gives a singular threshold without a
folded singularity].

Finally, if $\hat{\tau}_{\max}$ is positive and finite,
there may be `spurious' singular thresholds in phase portraits with a
folded singularity and $b\ge0$, or with no folded singularities, where
all trajectories starting on $S^a$ and near $F$ for $\hat{\tau}>0$
flow towards $F$. However, because $\hat{\tau}_{\max}$ is finite, some of these
trajectories will simply fail to reach $F$ by $\hat{\tau}_{\max}$.

It turns out that many examples of a singular threshold described
above, including the `spurious' singular threshold, can be eliminated
with a sensible assumption about $\lambda(\tau)$.

\subsubsection*{Part 2}
A more definitive statement about instability thresholds can be made
when $\tau_{\max}=\infty$, and the external input is
asymptotically constant, i.e. $\lambda(\tau)$ satisfies condition~(\ref{eq:asymptotic}).

Assume there is a singular threshold. On the one hand, it follows from
assumption (A1) and from condition~(\ref{eq:asymptotic}) that, for sufficiently large ${\tau}$,
trajectories started on $S^a$ and near $F$ must flow onto $S^a$ and
approach $\tilde{x}$.  On the other hand, a singular threshold
requires trajectories that start on $S^a$ and reach $F$. Hence, the
flow on $F$ in the desingularised system
(\ref{eq:fastd})--(\ref{eq:td}) must point towards $S^a$
for large values of $\tau$, and towards $S^r$ for lower values of
$\tau$. Such a change in the direction of the flow on $F$ requires a
folded singularity with $b>0$
in~(\ref{eq:fastn})--(\ref{eq:tn}). Hence, a folded saddle is
necessary for a singular threshold.

Assume there is a folded saddle singularity. There are two possible
situations. First, a folded saddle is the only folded
singularity. Second, a folded saddle is one of many folded
singularities.  In the second situation, assumption (A1)
and condition~(\ref{eq:asymptotic}) require that, typically, the
folded singularity with the largest $\tau$-component is a folded
saddle.  ``Typically'' excludes a folded saddle-node which is not
structurally stable. In both situations, there is a singular
threshold by the argument used for an isolated folded saddle in Part 1
of this proof. Hence, a folded saddle is sufficient for a singular
threshold.

\subsubsection*{Part 3}

In the last step of the proof we use theorems from canard theory
stating that the singular canards due to a folded
saddle~\cite[Th.~4.1]{Szmolyan2001}, a folded
node~\cite[Th.~4.1]{Szmolyan2001}\cite[Prop.~4.1]{Wechselberger2005},
and a folded saddle-node type I~\cite[Ths.~4.1 and~4.4]{Vo2014}, perturb
to {\it maximal canards} in~(\ref{eq:fast})--(\ref{eq:t}) with
$0<\delta \ll 1$.
Maximal canards are transverse, robust intersections of
two-dimensional attracting $S_{\delta}^a$ and repelling $S_{\delta}^r$
slow manifolds~\cite{Szmolyan2001,Wechselberger2005}. Such
intersections are possible in system~(\ref{eq:fast})--(\ref{eq:t})
because the slow manifolds $S_{\delta}^a$ and $S_{\delta}^r$ can be
extended across the fold~\cite{Desroches2012}. Starting on
$S_{\delta}^a$ and near the fold, trajectories jump off $S_{\delta}^a$
in the fast $x$-direction on one side of such intersections,
and flow
onto $S_{\delta}^a$ on the other side~\cite[Fig. 13]{Szmolyan2001}.
Thus, a singular threshold in system~(\ref{eq:fastp})--(\ref{eq:tp})
implies an instability threshold in
system~(\ref{eq:fast})--(\ref{eq:t}), and in the original
system~(\ref{eq:ofast})--(\ref{eq:oslow}).\\




\begin{thebibliography}{9}


\expandafter\ifx\csname urlstyle\endcsname\relax
  \providecommand{\doi}[1]{doi:\discretionary{}{}{}#1}\else
  \providecommand{\doi}{doi:\discretionary{}{}{}\begingroup
  \urlstyle{rm}\Url}\fi

\bibitem{Wieczorek2010} Wieczorek S, Ashwin P, Luke CM, Cox PM. 2010
Excitability in ramped systems: the compost-bomb instability.  {\em
Proc. R. Soc. A\/} {\bf 467}, 1243--1269.
(\doi{10.1098/rspa.2010.0485})

\bibitem{Lenton2008} Lenton T, Held H, Kriegler E, Hall J, Lucht W,
Rahmstorf S, Schellnhuber H.  2008 Tipping elements in the {E}arth's
climate system.  {\em PNAS\/} {\bf 105}, 1786--1793.
(\doi{10.1073/pnas.0705414105})

\bibitem{Stocker1997} Stocker TF, Schmittner A. 1997 Influence of
{CO2} emission rates on the stability of the thermohaline circulation.
{\em Nature\/} {\bf 388}, 862--865.

\bibitem{Leemans2004} Leemans R, Eickhout B. 2004 Another reason for
concern: regional and global impacts on ecosystems for different
levels of climate change.  {\em Global Envtl Change\/} {\bf 14},
219--228.  (\doi{10.1016/j.gloenvcha.2004.04.009})

\bibitem{Scheffer2008} Scheffer M, {van Nes} E, Holmgren M, Hughes
T. 2008 Pulse-driven loss of top-down control: the critical-rate
hypothesis.  {\em Ecosystems\/} {\bf 11}, 226--237.
(\doi{10.1007/s10021-007-9118-8})

\bibitem{Bridle2007} Bridle JR, Vines TH. 2007 Limits to evolution at
range margins: when and why does adaptation fail?  {\em Trends
Ecol. Evol.\/} {\bf 22}, 140--147.  (\doi{10.1016/j.tree.2006.11.002})

\bibitem{Izhikevich2007} Izhikevich E. 2007 {\em Dynamical Systems in
Neuroscience\/}.  Computational Neuroscience. MIT Press.

\bibitem{Hill1936} Hill AV. 1936 Excitation and accommodation in
nerve.  {\em Proc. R. Soc. Lond. B\/} {\bf 119}, 305--355.
(\doi{10.1098/rspb.1936.0012})

\bibitem{Biktashev2002} Biktashev VN. 2002 Dissipation of the
excitation wave fronts.  {\em Phys. Rev. Lett.\/} {\bf 89}, 168102.
(\doi{10.1103/PhysRevLett.89.168102})

\bibitem{Nene2012} Nene N, Garca-Ojalvo J, Zaikin A. 2012
Speed-dependent cellular decision making in nonequilibrium genetic circuits.  {\em PLoS ONE\/} {\bf 7}, 32779. 
(\doi{10.1371/journal.pone.0032779})

\bibitem{Ashwin2012} Ashwin P, Wieczorek S, Vitolo R, Cox PM. 2012
Tipping points in open systems: bifurcation, noise-induced and
rate-dependent examples in the climate system.  {\em
Phil. Trans. R. Soc. A\/} {\bf 370}, 1166--1184.
(\doi{10.1098/rsta.2011.0306})

\bibitem{Mitry2013} Mitry J, McCarthy M, Kopell N, Wechselberger
M. 2013 Excitable neurons, firing threshold manifolds and canards.
{\em J. Math. Neuro.\/} {\bf 3}, 12.  (\doi{10.1186/2190-8567-3-12})

\bibitem{Lobry1991} Beno{\^\i}t E. (ed) 1991 {\em Dynamic
Bifurcations.\/} Lecture Notes in Mathematics, vol. 1493.  Berlin,
Germany: Springer.

\bibitem{Luke2011} Luke CM, Cox PM. 2011 Soil carbon and climate
change: from the {J}enkinson effect to the compost{-}bomb instability.
{\em Eur. J. Soil Sci.\/} {\bf 62}, 5--12.
(\doi{10.1111/j.1365-2389.2010.01312.x})

\bibitem{Wechselberger2014} Wechselberger M, Mitry J, Rinzel J. 2013
Canard theory and excitability.  In {\em Nonautonomous Dynamical
Systems in the Life Sciences\/}, (eds PE~Kloeden, C~Poetzsche)
pp. 89--132. Lecture Notes in Mathematics, vol. 2102.  Springer
International Publishing.  (\doi{10.1007/978-3-319-03080-7_3})

\bibitem{Szmolyan2001} Szmolyan P, Wechselberger M. 2001 Canards in
{R}3.  {\em J. of Diff. Eqn.\/} {\bf 177}, 419--453.
(\doi{10.1006/jdeq.2001.4001})

\bibitem{Krupa2010} Krupa M, Wechselberger M. 2010 Local analysis near
a folded saddle-node singularity.  {\em J. Diff. Eqn.\/} {\bf 248},
2841--2888.  (\doi{10.1016/j.jde.2010.02.006})

\bibitem{Guckenheimer2008} Guckenheimer J. 2008 Return maps of folded
nodes and folded saddle-nodes.  {\em Chaos\/} {\bf 18}, 015108.
(\doi{10.1063/1.2790372})

\bibitem{Vo2014} Vo T, Wechselberger M.  Canards of folded saddle-node
type. Preprint submitted to {\em SIAM J. Math. Anal.}.

\bibitem{Fenichel1979} Fenichel N. 1979 Geometric singular
perturbation theory for ordinary differential equations.  {\em
J. Diff. Eqn.\/} {\bf 31}, 53--98.
(\doi{10.1016/0022-0396(79)90152-9})

\bibitem{Jones1995} Jones C. 1995 Geometric singular perturbation
theory.  In {\em Dynamical Systems\/}, (ed R~Johnson)
pp. 44--118. Lecture Notes in Mathematics, vol. 1609. Berlin, Germany:
Springer.  (\doi{10.1007/BFb0095239})

\bibitem{Roberts2013} Roberts A, Widiasih E, Jones CKRT, Wechselberger
M.  Mixed mode oscillations in a conceptual climate model. 2013.
ArXiv:1311.5182.

\bibitem{Cessi1994} Cessi P. 1994 A simple box model of stochastically
forced thermohaline flow.  {\em J. Phys. Oceanogr.\/} {\bf 24},
1911--1920.  (\doi{10.1175/1520-0485(1994)024<1911:ASBMOS>2.0.CO;2})

\bibitem{Pol1934} {van der Pol} B. 1934 The nonlinear theory of
electric oscillations.  {\em Proc. IRE\/} {\bf 22}, 1051--1086.
(\doi{10.1109/JRPROC.1934.226781})

\bibitem{Takens1976} Takens F. 1976 Constrained equations; a study of
implicit differential equations and their discontinuous solutions.  In
{\em Structural Stability, the Theory of Catastrophes, and
Applications in the Sciences\/}, (ed P~Hilton) pp. 143--234. Lecture
Notes in Mathematics, vol. 525. Berlin, Germany: Springer.
(\doi{10.1007/BFb0077850})

\bibitem{Dumortier1996} Dumortier F, Roussarie RH. 1996 {\em Canard
Cycles and Center Manifolds\/}.  {\em Mem. Amer. Math. Soc.\/}
No. 577.

\bibitem{Wechselberger2005} Wechselberger M. 2005 Existence and
bifurcations of canards in {R}3 in the case of a folded node.  {\em
SIAM J. App. Dy. Sys.\/} {\bf 4}, 101--139.  (\doi{10.1137/030601995})

\bibitem{Desroches2012} Desroches M, Guckenheimer J, Krauskopf B,
Kuehn C, Osinga HM, Wechselberger M.  2012 Mixed mode oscillations
with multiple timescales.  {\em SIAM Review\/} {\bf 54}, 211--288.
(\doi{10.1137/100791233})


\bibitem{Wechselberger2009} Wechselberger M, Weckesser W. 2009
Bifurcations of mixed-mode oscillations in a stellate cell model.
{\em Phys. D\/} {\bf 238}, 1598--1614.
(\doi{10.1016/j.physd.2009.04.017})

\bibitem{Desroches2010} Desroches M, Krauskopf B, Osinga HM.  2010 
Numerical continuation of canard orbits in slow-­fast dynamical systems.
{\em Nonlinearity\/} {\bf 23}, 739--765.
(\doi{10.1088/0951-7715/23/3/017})

\bibitem{Szmolyan2004} Szmolyan P, Wechselberger M. 2004
Relaxation oscillation in {R}3.
{\em J. Diff. Eqn.\/} {\bf 200}, 69--104.
(\doi{10.1016/j.jde.2003.09.010})


\end{thebibliography}
\end{document}